\documentclass[12pt,twoside]{amsart}

\usepackage[latin1]{inputenc}

\usepackage{amsmath,amsfonts,amssymb,amsthm,fullpage,bbm}

\usepackage{fancyhdr,graphicx,color}

\numberwithin{equation}{section}

\newtheorem{thm}{Theorem}[section]

\newtheorem{lem}[thm]{Lemma}

\newtheorem{cor}[thm]{Corollary}

\newtheorem{prop}[thm]{Proposition}

\newtheorem{rem}[thm]{Remark}

\newtheorem{dfn}[thm]{Definition}
\newcommand{\diam}{diam}
\newcommand{\aver}[1]{-\hskip-0.46cm\int_{#1}}
\newcommand{\avert}[1]{-\hskip-0.38cm\int_{#1}}
\newcommand{\ind}{1\hspace{-2.5 mm}{1}}
\makeatletter
\renewcommand\@biblabel[1]{#1.}
\makeatother

\begin{document}
\allowdisplaybreaks
\title{ Gagliardo-Nirenberg inequalities  on Manifolds}
\author{Nadine Badr}
\address{N.Badr
\\
Universit\'e de Paris-Sud, UMR du CNRS 8628
\\
91405 Orsay Cedex, France} \email{nadine.badr@math.u-psud.fr}

\date{\today}
\begin{abstract} We prove Gagliardo-Nirenberg inequalities on some classes of manifolds, Lie groups and graphs. 
\end{abstract}
\maketitle
\tableofcontents

\section{Introduction}
 Cohen-Meyer-Oru \cite{cohen1}, Cohen-Devore-Petrushev-Xu \cite{cohen2}, proved the following Gagliardo-Nirenberg type inequality 
\begin{equation}\label{C}
\|f\|_{1^*}\leq C\|\,|\nabla f|\,\|_{1}^{\frac{n-1}{n}}\|f\|_{B_{\infty,\infty}^{-(n-1)}}^{\frac{1}{n}}
\end{equation} 
for all $f\in W_{1}^{1}(\mathbb{R}^{n})$ ($1^*=\frac{n}{n-1}$). 
The proof of (\ref{C}) is involved and based on wavelet decompositions, weak type (1,1) estimates and interpolation results.

Using a simple method relying on weak type estimates and pseudo-Poincar\'{e} inequalities, Ledoux \cite{ledoux} obtained the following extension of (\ref{C}). He proved that for $1\leq p<l<\infty$ and for every $f\in W_{p}^{1}(\mathbb{R}^{n})$ 
\begin{equation}\label{led}
\|f\|_{l}\leq C \|\,|\nabla f|\,\|_{p}^{\theta} \|f\|_{B_{\infty,\infty}^{\frac{\theta}{\theta-1}}}^{1-\theta}
\end{equation}
 where $\theta= \frac{p}{l}$ and $C>0$ only depends on $l,\,p $ and $n$.
 
 In the same paper, he extended (\ref{led}) to the case of Riemannian manifolds. If $p=2$ he observed that (\ref{led}) holds without any assumption on $M$. If $p\neq 2$ he assumed that the Ricci curvature is non-negative and obtained (\ref{led}) with $C>0$ only depending on $l,\,p$ when $1\leq p\leq 2$ and on $l,\,p $ and $n$  when $2<p<\infty$.
 
 He also proved that a similar inequality holds on $\mathbb{R}^{n}$, Riemannian manifolds with non-negative Ricci curvature, Lie groups and Cayley graphs, replacing the $B_{\infty,\infty}^{\frac{\theta}{\theta-1}}$ norm by the $M_{\infty}^{\frac{\theta}{\theta-1}}$ norm (see definitions below). 

Note that these two versions of Gagliardo-Nirenberg inequalities extend the classical Sobolev inequality in the Euclidean case:
\begin{equation}\label{SI}
\|f\|_{p^*}\leq C\|\,|\nabla f|\,\|_{p}
\end{equation}
with $\frac{1}{p^*}=\frac{1}{p}-\frac{1}{n}$  holds on $\mathbb{R}^{n}$ for every $f\in W_{p}^{1}(\mathbb{R}^{n})$ and for every $1\leq p<n$.
 \\
 In the Riemannian case it is not generally true that (\ref{led}) or (\ref{C}) imply (\ref{SI}), without
additional assumptions on the manifold (cf. Proposition \ref{SIP} below). On the other
hand we will now show examples of Riemannian manifolds where (\ref{SI}) holds independently
of (\ref{led}). It is clear that (\ref{SI}) holds on a compact Riemannian $n$-manifold $M$. 
As an example of complete non-compact Riemannian manifold satisfying (\ref{SI}), we can consider a complete Riemannian $n$-manifold $M$ with non-negative Ricci curvature. If there exists $v>0$ such that for all $x\in M$, $\mu(B(x,1))\geq v$, then $M$ satisfies (\ref{SI}). Here $\mu(B(x,1))$ is the Riemannian volume of the open ball $B(x,1)$. For more general cases where we have (\ref{SI}) for some $p$'s depending on the hypotheses, see \cite{saloff3}. 
 Note that if $(\ref{SI})$ holds for some $1\leq p<n$, then it holds for all $p\leq q<n$ (see \cite{saloff3}, Chapter 3).

We have also non-linear versions of Gagliardo-Nirenberg inequalities proved by Rivi\`{e}re-Strzelecki \cite{riviere1}, \cite{strzelecki}. They got for every $f\in C_{0}^{\infty}(\mathbb{R}^{n})$
\begin{equation}\label{riv}
\int_{\mathbb{R}^{n}}|\nabla f|^{p+2}\leq C\|f\|_{BMO}^{2}\int_{\mathbb{R}^{n}}|\nabla^{2}f|^{2}|\nabla f|^{p-2}.
\end{equation}

They applied this inequality and obtained a regularity property for solutions of non-linear elliptic equations of type
$$
-div (|\nabla u|^{p-2}\nabla u)=G(x,u,\nabla u)
$$
where $G$ grows as $|\nabla u|^{p}$.

Recently, Martin-Milman \cite{martin} developed a new symmetrization approach to obtain the Gagliardo-Nirenberg inequalities (\ref{led}) and, therefore the Sobolev inequalities (\ref{SI}) in $\mathbb{R}^{n}$. They also proved a variant of (\ref{riv}). The method of \cite{martin} to prove (\ref{led}) is different from that of Ledoux. It relies essentially on an interpolation result for Sobolev spaces and pseudo-Poincar\'e inequalities in the Euclidean case.
 
 In this paper, we prove analogous results on Riemannian manifolds, Lie groups and graphs making some additional hypotheses on these spaces. For this purpose, we will adapt Martin and Milman's method and make use of our interpolation results in \cite{badr1}. More precisely we obtain in the case of Riemannian manifolds:
\begin{thm}\label{G} Let $M$ be a complete non-compact Riemannian manifold satisfying $(D)$ and $(P_{q})$  for some $1\leq q<\infty$. Moreover, assume that $M$ satisfies the pseudo-Poincar\'{e} inequalities $(P'_{q})$ and $(P'_{\infty})$. Consider $\alpha<0$. Then, there exists $C>0$ such that for every $f \in (W_{q}^{1}+W_{\infty}^{1})\cap B_{\infty,\infty}^{\alpha}$  with $f^{*}(\infty)=0$ and $|\nabla f|^{*}(\infty)=0$, we have
\begin{equation}\label{B}
|f|^{q**\frac{1}{q}}(s)\leq C|\nabla f|^{q**\frac{|\alpha|}{q(1+|\alpha|)}}(s)\|f\|_{B_{\infty,\infty}^{\alpha}}^{\frac{1}{1+|\alpha|}}.
\end{equation}
\end{thm}
Above and from now on,  $|f|^{q**\frac{1}{q}}$ means $(|f|^{q**})^{\frac{1}{q}}$.
Recall that for every $t>0$
$$
f^{*}(t)=\inf\left\lbrace\lambda;\mu(\{|f|>\lambda\right\rbrace)\leq t\};
$$

$$
f^{*}(\infty)=\inf\left\lbrace\lambda;\mu(\{|f|>\lambda\right\rbrace)<\infty\}
$$
and
$$
f^{**}(t)=\frac{1}{t}\int_{0}^{t}f^{*}(s)ds.
$$
Using this symmetrization result we prove
\begin{thm} \label{L} Let $M$ be a complete Riemannian manifold satisfying the hypotheses of Theorem \ref{G}. Then (\ref{led}) holds for all $q\leq p<l<\infty$.
\end{thm}
\begin{cor}\label{R} Let $M$ be a Riemannian manifold with non-negative Ricci curvature. Then (\ref{led}) holds for all $1\leq p<l<\infty$.
\end{cor} 
This corollary is exactly what Ledoux proved \cite{ledoux}. We obtain further generalizations: 

\begin{cor}\label{lp} Consider a complete Riemannian manifold $M$ satisfying $(D)$, $(P_{1})$ and assume that there exists $C>0$ such that for every $x,\,y\in M$ and $t>0$
\begin{equation*}\tag{$G$}
|\nabla_{x}p_{t}(x,y)|\leq \frac{C}{\sqrt{t}\mu(B(y,\sqrt{t}))}.
\end{equation*}
Then inequality (\ref{led}) holds for all $1\leq p<l<\infty$.
\end{cor}
Note that a Lie group of polynomial growth satisfies the hypotheses of Corollary \ref{lp} (see \cite{coulhon7}). Hence it verifies (\ref{led}) for all $1\leq p<l<\infty$. 
\\
 Another example of a space satisfying the hypotheses of Corollary \ref{lp}
  is given by taking a Galois covering manifold of a compact manifold whose deck transformation group has polynomial growth (see \cite{dungey2}). We can also take the example of a Cayley graph of a finitely generated group (see \cite{coulhon6}, \cite{saloff3}).\\
 \\
We also get the following Corollary:
\begin{cor}\label{P2} Let $M$ be a complete Riemannian manifold satisfying $(D)$ and $(P_{2})$. Then (\ref{led}) holds for all $2\leq p<l<\infty$.
\end{cor}
Note that $(P'_{2})$ is always satisfied. Hence, by Ledoux's method, inequality (\ref{led}) with $p=2$ needs no assumption on $M$ (see \cite{ledoux}). So our results are only interesting when $p\neq 2$.

\paragraph{\textbf{Local version:}} Let $M$ be a complete Riemannian manifold satisfying a local doubling property $(D_{loc})$ and a local Poincar\'e inequality $(P_{qloc})$  --we restrict our definitions to small balls--. Moreover assume that $M$ admits a local version of pseudo-Poincar\'e inequalities $(P'_{qloc})$, $(P'_{\infty loc})$: by $(P'_{rloc})$ we mean
$$
\|f-e^{t\Delta}f\|_{r}\leq Ct^{\frac{1}{2}}\left(\|f\|_{r}+\|\,|\nabla f|\,\|_{r}\right).
$$
In this context, the following local version of (\ref{led}) holds: for every $q\leq p<l<\infty$ and $f\in W_{p}^{1}$
\begin{equation}\label{led"}
\|f\|_{l}\leq C \left(\|f\|_{p}+\|\,|\nabla f|\,\|_{p}\right)^{\theta} \|f\|_{B_{\infty,\infty}^{\frac{\theta}{\theta-1}}}^{1-\theta}.
\end{equation}
\\
 
In the following theorem, we show a variant of Theorem \ref{G} replacing the Besov norm by the Morrey norm. In the Euclidean case, the Morrey space is strictly smaller than the Besov space. Therefore, the following Theorem \ref{M} (resp. Corollary \ref{L'}) is weaker than Theorem \ref{G} (resp. Theorem \ref{L}). In contrast, on Riemannian manifolds, the Besov and Morrey spaces are not comparable in general.
\begin{thm}\label{M} Let $M$ be a complete non-compact Riemannian manifold satisfying $(D)$ and $(P_{q})$ for some $1\leq q<\infty$. Consider $q\leq p<\infty$ and $\alpha<0$. Then, for every $f\in (W_{q}^{1}+W_{\infty}^{1})\cap M_{\infty}^{\alpha}$ we have
$$
|f|^{q**\frac{1}{q}}(s)\leq C|\nabla f|^{q**\frac{|\alpha|}{q(1+|\alpha|)}}(s)\|f\|_{M_{\infty}^{\alpha}}^{\frac{1}{1+|\alpha|}}.
$$
\end{thm}
\begin{cor}\label{L'} Under the hypotheses of Theorem \ref{M}, let $q_{0}= \inf\left\lbrace q \in [1,\infty[: (P_{q}) \textrm{ holds }\right\rbrace$ and consider $q_{0}<p<l<\infty$\footnote{if $q_{0}=1$, we allow $ 1\leq p<l<\infty$}. 
 Then, for every $f\in W_{p}^{1}$, we have
\begin{equation}\label{led'}
\|f\|_{l}\leq C \|\,|\nabla f|\,\|_{p}^{\theta} \|f\|_{M_{\infty}^{\frac{\theta}{\theta-1}}}^{1-\theta}.
\end{equation}
\end{cor} 

\noindent Ledoux \cite{ledoux} showed that (\ref{led'}) holds on any unimodular Lie group equipped with a left invariant Riemannian metric and the associated Haar measure. Once again, this is due to the fact that his method uses essentially the pseudo-Poincar\'e inequalities $(P^{\prime\prime}_{p})$, which hold on such a group for all $1\leq p\leq \infty$ (see \cite{saloff3}). With our method, we only get the local version of (\ref{led'}), namely the analog of (\ref{led"}). However notice that we prove (\ref{led'}) in its full strength for Lie groups of polynomial growth.\\

Let us compare our result with Ledoux's one. Our hypotheses are stronger, we assume in addition of the pseudo-Poincar\'e inequality --which is the only assumption of Ledoux-- $(D)$ and $(P_q)$ but recover most of his examples. Moreover we obtain Corollary \ref{lp} which gives us more examples as we have seen in the introduction. For instance, on Lie groups, Ledoux only mentioned in his paper the Morrey version while Corollary \ref{lp} yield (\ref{led}) on Lie groups with polynomial growth for every $1\leq p<l<\infty$. We get also the interpolation of his inequality (\ref{led}). Since it is not known if the pseudo-Poincar\'e inequalities interpolate, his method gives (\ref{led}) (resp. (\ref{led'})) for the same exponent $p$ of pseudo-Poincar\'e inequality. With our method, we get (\ref{led}) (resp. (\ref{led'})) for every $p\geq q$.
\\
 
We finish with the following non-linear Gagliardo-Nirenberg theorem:
\begin{thm}\label{S} Let $M$ be a complete non-compact Riemannian manifold satisfying $(D)$ and $(P_{q})$ for some $1\leq q< \infty$. Moreover, assume that $M$ satisfies $(P'_{q})$ and $(P'_{\infty})$. Let $p\geq \max(2,q)$. Then for every $f\in C^{\infty}_{0}(M)$
 $$
 \int_{M}|\nabla f|^{p+1}d\mu\leq C \|f\|_{B_{\infty,\infty}^{-1}}\int_{M}|\nabla^{2}f|^{2}|\nabla f|^{p-2}d\mu.
 $$
\end{thm}

The paper is organized as follows. In section 2, we give the definitions on a Riemannian manifold of Besov and Morrey spaces, Sobolev spaces, doubling property, Poincar\'{e} and pseudo-Poincar\'{e} inequalities.
In section 3, we show how to obtain under our hypotheses Ledoux's inequality (\ref{led}) and different Sobolev inequalities. Section 4 is devoted to prove Theorem \ref{G} and Theorem \ref{M}. In section 5 we give another symmetrization inequality. Finally we prove Theorem \ref{S} in section 6.
\\

\textit{Acknowledgements.} I would like to thank my Ph.D advisor P. Auscher for his comments and advice about the topic of this paper. I am also indebted to J. Mart\'{i}n and M. Milman for the useful discussions I had with them, especially concerning Theorem \ref{EK'H}.

\section{Preliminaries} Throughout this paper $C$ will be a constant
that may change from an inequality to another and we will use $u\sim
v$ to say that there exist two constants $C_{1}$,$C_{2}>0$ such that $C_{1}u\leq v\leq
C_{2}u$.

 Let $M$ be a complete non-compact Riemannian manifold. We write $\mu$ for the Riemannian measure on $M$, $\nabla$ for the Riemannian gradient, $|\cdot|$ for the length on the tangent space (forgetting the subscript $x$ for simplicity) and $\|\cdot\|_{p}$ for the norm on $ L_{p}(M,\mu)$, $1 \leq p\leq +\infty$. Let $P_{t}=e^{t\Delta}$, $t\geq 0$, be the heat semigroup on $M$ and $p_{t}$ the heat kernel.
\subsection{Besov and Morrey spaces}
For $\alpha<0$, we introduce the Besov norm
$$
\|f\|_{B_{\infty,\infty}^{\alpha}}=\sup\limits_{t>0}t^{-\frac{\alpha}{2}}\|P_{t}f\|_{\infty}<\infty
$$
for measurable functions $f$ such that this makes sense and say $f\in B_{\infty,\infty}^{\alpha}$ (we shall not try here to give the most general definition of the Besov space).
\begin{lem} We have  for every $f\in B_{\infty,\infty}^{\alpha}$
\begin{equation}\label{eB}
\|f\|_{B_{\infty,\infty}^{\alpha}}\sim \sup_{t>0}t^{-\frac{\alpha}{2}}\|P_{t}(f-P_{t}f)\|_{\infty}.
\end{equation}
\end{lem}
\begin{proof}
It is clear that $\sup_{t>0}t^{-\frac{\alpha}{2}}\|P_{t}(f-P_{t}f)\|_{\infty}\leq (1+2^{\frac{\alpha}{2}})\|f\|_{B_{\infty,\infty}^{\alpha}}$. On the other hand
$$
 t^{-\frac{\alpha}{2}}P_{t}f=t^{-\frac{\alpha}{2}}(P_{t}f-P_{2t}f)+2^{\frac{\alpha}{2}}\,(2t)^{-\frac{\alpha}{2}}P_{2t}f.
$$
By taking the supremun over all $t>0$, we get
$$
\|f\|_{B_{\infty,\infty}^{\alpha}}\leq \sup_{t>0}t^{-\frac{\alpha}{2}}\|P_{t}(f-P_{t}f)\|_{\infty}+2^{\frac{\alpha}{2}}\|f\|_{B_{\infty,\infty}^{\alpha}}.
$$
Thus, $\|f\|_{B_{\infty,\infty}^{\alpha}}\leq \frac{1}{1-2^{\frac{\alpha}{2}}}\sup_{t>0}t^{-\frac{\alpha}{2}}\|P_{t}(f-P_{t}f)\|_{\infty}$.
\end{proof}
For $\alpha<0$, the Morrey space $M_{\infty}^{\alpha}$ is the space of locally integrable functions $f$ for which the Morrey norm  
$$
\|f\|_{M_{\infty}^{\alpha}}:=\sup\limits_{r>0,\,x\in M}r^{-\alpha}|f_{B(x,r)}|<\infty
$$
where $f_{B}:=\avert{B}fd\mu=\frac{1}{\mu(B)}\int_{B}fd\mu.$
 \subsection{ Sobolev spaces on Riemannian manifolds}
 \begin{dfn}[\cite{aubin1}] Let $M$ be a $C^{\infty}$ Riemannian manifold
 of dimension $n$. Write $E^{1}_{p}$ for the vector space of $C^{\infty}$ functions $\varphi
$ such that $\varphi $ and $|\nabla\varphi|\in L_{p},
\,1\leq p< \infty$. We define the non-homogeneous Sobolev space $W^{1}_{p}$ as the completion of  $E^{1}_{p}$ for the norm
$$
\|\varphi\|_{W^{1}_{p}}=\|\varphi\|_{p}+\|\,|\nabla\varphi|\,\|_{p}.
$$
We denote $W^{1}_{\infty}$ for the set of all bounded Lipschitz functions on $M$.
\end{dfn}
\begin{prop}\label{CDW} (\cite{aubin1}) Let $M$ be a complete Riemannian manifold. Then $ C^{\infty}_{0}$ is dense in $W^{1}_{p}$ for $1\leq p<\infty$.
\end{prop}
\begin{dfn}Let $M$ be a $C^{\infty}$ Riemannian manifold of dimension $n$.
For $1\leq p\leq \infty$, we define $\overset{.}{E_{p}^{1}}$ to be the vector space of distributions $ \varphi $
 with $|\nabla \varphi |\in L_{p}$, where $\nabla \varphi$ is the distributional gradient of $ \varphi$. It is well known that the elements of $\overset{.}{E_{p}^{1}}$ are in  $L_{p,loc} $. We equip $\overset{.}{E_{p}^{1}}$  with the semi norm 
 $$ 
 \|\varphi\|_{\overset{.}{E_{p}^{1}}}=\|\,|\nabla \varphi|\,\|_{p}.
 $$
 \end{dfn}
 \begin{dfn} We define the homogeneous Sobolev space $\overset{.}{W_{p}^{1}}$ as the quotient space $\overset{.}{E_{p}^{1}}/\mathbb{R}$.
\end{dfn}
\begin{rem} For all $\varphi \in \overset{.}{E_{p}^{1}}$, \, $\|\overline{\varphi}\|_{\overset{.}{W_{p}^{1}}}=\|\,|\nabla \varphi|\,\|_{p}$.
\end{rem}
 \subsection{Doubling property and Poincar\'{e} inequalities} 
\begin{dfn}[Doubling property] Let $(M,d,\mu)$ be a Riemannian manifold. Denote by $B(x, r)$ the open ball of center $x\in M$ and radius $r>0$. One says that $M$ satisfies the doubling property $(D)$ if there exists a constant $C_{d}>0$ such that for all $x\in M,\, r>0$ we have
\begin{equation*}\tag{$D$}
\mu(B(x,2r))\leq C_{d} \mu(B(x,r)).
\end{equation*}
\end{dfn}
 Observe that if $M$ satisfies $(D)$ then
$$
 \diam(M)<\infty\Leftrightarrow\,\mu(M)<\infty\; \textrm{(see \cite{ambrosio1})}.
  $$
 
\begin{dfn}[Poincar\'{e} inequality] A complete Riemannian manifold $M$ admits a Poincar\'{e} inequality $(P_{q})$ for some $1\leq q<\infty$ if there exists a constant $C$ such that for all $f\in C^{\infty}_{0}$  and for every ball $B$ of $M$ of radius $r>0$, we have
\begin{equation*}\tag{$P_{q}$}
\Bigl( \avert{B}|f-f_{B}|^{q}d\mu\Bigr)^{\frac{1}{q}}\leq Cr \Bigl(\aver{B}|\nabla f|^{q}d\mu\Bigr)^{\frac{1}{q}}
\end{equation*} 
\end{dfn}
\begin{rem} Since $C^{\infty}_{0}$ is dense in $W_{q}^{1}$, if $M$ admits $(P_{q})$ for all $f\in C^{\infty}
_{0}$ then $(P_{q})$ holds for all  $f \in W_{q}^{1}$. In fact, by Theorem 1.3.4 in \cite{keith1}, $M$ admits $(P_{q})$ for all $f\in \dot{E}_{q}^{1}$.
\end{rem}
 The following recent result from Keith-Zhong \cite{keith3} improves the exponent of Poincar\'e inequality:
\begin{thm}\label{kz} Let $(X,d,\mu)$ be a complete metric-measure space with $\mu$ locally doubling
and admitting a local Poincar\'{e} inequality $(P_{q})$, for  some $1< q<\infty$.
Then there exists $\epsilon >0$ such that $(X,d,\mu)$ admits
$(P_{p})$ for every $p>q-\epsilon$.
\end{thm} 
\begin{dfn}[Pseudo-Poincar\'{e} inequality for the heat semigroup]
A Riemannian manifold $M$ admits a pseudo-Poincar\'{e} inequality  for the heat semigroup $(P'_{q})$ for some $1\leq q<\infty$ if there exists a constant $C$ such that for all $f\in C^{\infty}_{0}$ and all $t>0$, we have
\begin{equation*}\tag{$P'_{q}$}
\|f-P_{t}f\|_{q}\leq Ct^{\frac{1}{2}}\|\,|\nabla f|\,\|_{q}.
\end{equation*} 
$M$ admits a pseudo-Poincar\'{e} inequality $(P'_{\infty})$ if there exists $C>0$ such that for every bounded Lipschitz function $f$ we have
\begin{equation*}\tag{$P'_{\infty}$}
\|f-P_{t}f\|_{\infty}\leq Ct^{\frac{1}{2}}\|\,|\nabla f|\,\|_{\infty}.
\end{equation*} 
\end{dfn}
\begin{rem}  Again by density of  $C^{\infty}_{0}$  in $W_{q}^{1}$, if $M$ admits $(P'_{q})$ for some $1\leq q<\infty$ for all $f\in C^{\infty}
_{0}$ then $M$ admits $(P'_{q})$ for all  $f \in W_{q}^{1}$.
\end{rem}
\begin{dfn}[Pseudo-Poincar\'{e} inequality for averages]
A complete Riemannian manifold $M$ admits a pseudo-Poincar\'{e} inequality  for averages $(P^{\prime\prime}_{q})$ for some $1\leq q<\infty$ if there exists a constant $C$ such that for all $f\in C^{\infty}_{0}$ and for every ball $B$ of $M$ of radius $r>0$, we have
\begin{equation*}\tag{$P^{\prime\prime}_{q}$}
\|f-f_{B(.,r)}\|_{q}\leq Cr\|\,|\nabla f|\,\|_{q}.
\end{equation*} 
\end{dfn}
\begin{rem}(Lemma 5.3.2 in \cite{saloff3}) If $M$ is a complete Riemannian manifold satisfying $(D)$ and $(P_{q})$ for some $1\leq q<\infty$, then it satisfies $(P^{\prime\prime}_{q})$. Hence $(P^{\prime\prime}_{q})$ holds for all $f\in\overset{.}{E_{q}^{1}}$.
\end{rem}
  
\section{Ledoux's and Sobolev inequalities}
\subsection{Ledoux's inequality}
 \begin{proof}[Proof of Theorem \ref{L}] Let us show how to obtain Ledoux's inequality (\ref{led}) from Theorem \ref{G}.
 Consider $M$ satisfying the hypotheses of Theorem \ref{G} and take $q<p<l$. From (\ref{B}), we see that 
$$
\|\,|f|^{q**\frac{1}{q}}\|_{X}\leq C \|\,|\nabla f|^{q**\frac{1}{q}}\|_{X_{\frac{|\alpha|}{1+|\alpha|}}}^{\frac{|\alpha|}{1+|\alpha|}} \|f\|_{B_{\infty,\infty}^{\alpha}}^{\frac{1}{1+|\alpha|}}
$$ 
with $X=L_{l}$  which is a rearrangement invariant space (see \cite{bennett}, section 2 of \cite{martin}) and
$$
X_{a}=\left\lbrace f:|f|^{a}\in X, \textrm {with } \|f\|_{X_{a}}=\|\,|f|^{a}\,\|_{X}^{\frac{1}{a}}\right\rbrace.
$$ 
By taking $\alpha=\frac{p}{p-l}$ we obtain (\ref{led}) for $p>q$. For $q=p$, note that (\ref{B}) implies the weak type inequality $(q,l)$, which is $\mu(\{|f|>\lambda\})\leq \left(\frac{C}{\lambda}\|\,|\nabla f|\,\|_{q}^{\frac{q}{l}}\|f\|_{B_{\infty,\infty}^{\alpha}}^{1-\frac{q}{l}}\right)^{l}$. Consequently the strong type $(q,l)$, which is $\|f\|_{l}\leq C\|\,|\nabla f|\,\|_{q}^{\frac{q}{l}}\|f\|_{B_{\infty,\infty}^{\alpha}}^{1-\frac{q}{l}}$, follows by Maz'ya's truncation principle (see \cite{hajlasz9}, \cite{ledoux}). 
\end{proof}
\begin{proof}[Proof of Corollary \ref{R}]
Remark that Riemannian manifolds with non-negative Ricci curvature satisfy $(D)$ (with $C_{d}=2^n$) ,$\,(P_{1})$. They also satisfy  $(P'_{p})$ for all $1\leq p\leq \infty$, where the constant $C$ is numerical for $1\leq p\leq 2$ and only depends on $n$  for $2<p\leq \infty$ (see \cite{ledoux}). Thus Theorem \ref{L} applies on such manifolds with $q=1$.
\end{proof} 
Before we prove Corollary \ref{lp}, we give the following two lemmas. Let $2<p\leq \infty$. Consider the  following condition: there exists $C>0$ such that for every $t>0$
\begin{equation}\tag{$G_{p}$}
\|\,|\nabla e^{t\Delta}|\,\|_{p\rightarrow p}\leq \frac{C}{\sqrt{t}}.
\end{equation}
\begin{lem}\label{l1}(\cite{dungey1}) Let $M$ be a complete Riemannian manifold $M$ satisfying $(D)$ and the Gaussian heat kernel upper bound, that is, there exist $C,c>0$ such that for every $x,\,y\in M$ and $t>0$ 
\begin{equation}\label{gp}
p_{t}(x,y)\leq \frac{C}{\mu(B(y,\sqrt{t}))}e^{-c\frac{d^{2}(x,y)}{t}}.
\end{equation}
Then $(G)$  holds if and only if $(G_{\infty})$ holds.
\end{lem}
\begin{lem}\label{l2} Let $M$ be a complete Riemannian manifold. If the condition $(G_{p})$ holds for some $1< p\leq \infty$ then $M$ admits a pseudo-Poincar\'{e} inequality $(P'_{p'})$, $p'$ being the conjugate of $p$ $(\frac{1}{p}+\frac{1}{p'}=1)$.
\end{lem}
\begin{proof}
For $f\in C_{0}^{\infty}$, we have 
$$ f-e^{t\Delta}f=-\int_{0}^{t}\Delta e^{s\Delta}f\,ds.
$$
Remark that $(G_{p})$ gives us that $\|\Delta e^{s\Delta}f \|_{p'}\leq \frac{C}{\sqrt{s}}\|\,|\nabla f|\,\|_{p'}$. Indeed
\begin{align*}
\|\Delta e^{s\Delta}f\|_{p'}&=\sup_{\|g\|_{p}=1}\int_{M}\Delta e^{s\Delta}f \,g\,d\mu
\\
&=\sup_{\|g\|_{p}=1}\int_{M}f \,\Delta e^{s\Delta}g\,d\mu
\\
&=\sup_{\|g\|_{p}=1}\int_{M}\nabla f.\nabla e^{s\Delta}g\,d\mu
\\
&\leq \|\,|\nabla f|\,\|_{p'}\,\sup_{\|g\|_{p}=1}\|\,|\nabla e^{s\Delta}g|\,\|_{p}
\\
&\leq \frac{C}{\sqrt{s}}\|\,|\nabla f|\,\|_{p'}.
\end{align*}
Therefore 
$$\|f-e^{t\Delta}f\|_{p'}\leq C\|\,|\nabla f|\,\|_{p'}\int_{0}^{t}\frac{1}{\sqrt{s}}ds =C\sqrt{t}\|\,|\nabla f|\,\|_{p'}
$$
which finishes the proof of the lemma.
\end{proof}
\begin{proof}[Proof of Corollary \ref{lp}] The fact that $M$ satisfies $(D)$ and admits $(P_{1})$, hence $(P_{2})$, gives the Gaussian heat kernel upper bound (\ref{gp}). Since $(G)$ holds, Lemma \ref{l1} asserts that $(G_{\infty})$ holds too. Applying Lemma \ref{l2} it comes that $M$ admits a pseudo-Poincar\'{e} inequality $(P'_{1})$. We claim that $(P'_{\infty})$ holds on $M$. Indeed, (\ref{gp}) yields
\begin{align*}
\|f-e^{t\Delta }f\|_{\infty}&\leq \sup_{x\in M}\int_{M}|f(x)-f(y)|p_{t}(x,y)d\mu(y)
\\
&\leq C\|\,|\nabla f|\,\|_{\infty} \sup_{x\in M}\frac{1}{\mu(B(x,\sqrt{t}))}\int_{M}d(x,y)e^{-c\frac{d^{2}(x,y)}{t}}d\mu(y)
\\
&\leq C\sqrt{t}\|\,|\nabla f|\,\|_{\infty}\sup_{x\in M}\frac{1}{\mu(B(x,\sqrt{t}))}\int_{M}e^{-c'\frac{d^{2}(x,y)}{t}}d\mu(y)
\\
&\leq C\sqrt{t}\|\,|\nabla f|\,\|_{\infty}\sup_{x\in M}\frac{1}{\mu(B(x,\sqrt{t}))}\mu(B(x,\sqrt{t}))
\\
&= C\sqrt{t}\|\,|\nabla f|\,\|_{\infty}
\end{align*}
where the last estimate is a straightforward consequence of $(D)$. Therefore, we have all we need to apply Theorem \ref{G} with $q=1$. The inequality (\ref{led}) for all $1\leq p<l<\infty$ follows then by Theorem \ref{L}.
\end{proof}
\begin{rem} Under the hypotheses of Corollary \ref{lp}, Theorem \ref{M} and Theorem \ref{S} also hold.
\end{rem}
 \begin{proof}[Proof of Corollary \ref{P2}] First we know that $(G_{2})$ always holds on $M$. Then $(P'_{2})$ holds by Lemma \ref{l2}. Secondly $(D)$ and $(P_{2})$ yields $(P'_{\infty})$ as we have just seen above. Hence Theorem \ref{L} applies with $q=2$. 
\end{proof}
\subsection{The classical Sobolev inequality}
\begin{prop}\label{SIP} Consider a complete non-compact Riemannian manifold satisfying the hypotheses of Theorem \ref{G} and assume that $1\leq q<\nu$ with $\nu>0$. From (\ref{led}) and under the heat kernel bound $\|P_{t}\|_{q \rightarrow \infty}\leq Ct^{-\frac{\nu}{2q}}$, one recovers the classical Sobolev inequality
\begin{equation*}
\|f\|_{q^*}\leq C \|\,|\nabla f|\,\|_{q}
\end{equation*}
with $\frac{1}{q^*}=\frac{1}{q}-\frac{1}{\nu}$.
Consequently, we get
\begin{equation*}
\|f\|_{p^*}\leq C \|\,|\nabla f|\,\|_{p}
\end{equation*}
with $\frac{1}{p^*}=\frac{1}{p}-\frac{1}{\nu}$ for $q\leq p<\nu$.
\end{prop}
\begin{proof} Recall that $\|f\|_{B_{\infty,\infty}^{\alpha}}\sim \sup_{t>0}t^{-\frac{\alpha}{2}}\|P_{t}(f-P_{t}f)\|_{\infty}$.
 The pseudo-Poincar\'{e} inequality $(P'_{q})$, (\ref{led}) and the heat kernel bound $\|P_{t}\|_{q\rightarrow \infty}\leq Ct^{-\frac{\nu}{2q}}$ yield
$$
\|f\|_{q^*}\leq C \|\,|\nabla f|\,\|_{q}^{\theta} \left(\sup_{t>0} t^{-\frac{1}{2}}\|f-P_{t}f\|_{q}\right)^{1-\theta}
\leq C \|\,|\nabla f|\,\|_{q}.
$$
Thus we get (\ref{SI}) with $p=q<\nu$ and $\frac{1}{q^*}=\frac{1}{q}-\frac{1}{\nu}$.
\end{proof}
\subsection{Sobolev inequalities for Lorentz spaces}\label{SL}
 For $1\leq p\leq \infty$, $0\leq r<\infty$ we note $L(p,r)$ the Lorentz space of functions $f$ such that
 $$
 \|f\|_{L(p,r)}=\left(\int_{0}^{\infty}(f^{**}(t)t^{\frac{1}{p}})^{r}\frac{dt}{t}\right)^{\frac{1}{r}}<\infty
 $$
 and 
 $$
 \|f\|_{L(p,\infty)}=\sup_{t}t^{\frac{1}{p}}f^{*}(t)<\infty.
 $$

Consider a complete non-compact Riemannian manifold $M$ satisfying $(D)$ 
and $(P_{q})$ for some $1\leq q<\infty$. Moreover, assume that the following global growth condition 
\begin{equation}\label{GM}
\mu(B)\geq Cr^{\sigma}
\end{equation}
holds for every ball $B\subset M$ of radius $r>0$ and for some $\sigma>q$ (Remark that $\sigma\geq n$). Using Remark 4 in \cite{kalis}, we get
\begin{equation}\label{K}
f^{**}(t)-f^{*}(t)\leq Ct^{\frac{1}{\sigma}}|\nabla f|^{q**\frac{1}{q}}(t)
 \end{equation}
 for every $f\in \overset{.}{E_{q}^{1}}$. 
  We can write (\ref{K}) as
 \begin{equation}\label{k}
 f^{**}(t)-f^{*}(t)\leq \left[Ct^{\frac{1}{\sigma}}|\nabla f|^{q**\frac{1}{q}}(t)\right]^{1-\theta}(f^{**}(t)-f^{*}(t))^{\theta},\; 0\leq \theta\leq 1.
 \end{equation}
 Take $\frac{1}{r}=\frac{1-\theta}{p*}+\frac{\theta}{l},\,\frac{1}{m}=\frac{1-\theta}{m_{0}}+\frac{\theta}{m_{1}}$ with $ 0\leq \theta\leq 1$, $\sigma\geq p>q$, $m_{0}\geq q$ and $\frac{1}{p^{*}}=\frac{1}{p}-\frac{1}{\sigma}$. Then from (\ref{k}) and H\"{o}lder's inequality, we obtain the following Gagliardo-Nirenberg inequality for Lorentz spaces
 \begin{equation}\label{SLI}
 \|f\|_{L(r,m)}\leq C \|\,|\nabla f|\,\|_{L(p,m_{0})}^{1-\theta}\|f\|_{L(l,m_{1})}^{\theta}.
 \end{equation}
 We used also the fact that for $1<p\leq \infty$ and $1\leq r\leq \infty$
 $$\|f\|_{L(p,r)}\sim \left[\int_{0}^{\infty} \left(t^{\frac{1}{p}}f^{*}(t)\right)^{r}\frac{dt}{t}\right]^{\frac{1}{r}}
 $$ 
to obtain the term $\|\,|\nabla f|\,\|_{L(p,m_{0})}$ (see \cite{stein3}, Chapter 5, Theorem 3.21).\\
If we take $\theta=0$ and $m_{0}=m=p$, $r=p^{*}$, (\ref {SLI}) becomes
\begin{equation}\label{SLI'}
 \|f\|_{L(p^{*},p)}\leq C \|\,|\nabla f|\,\|_{p}.
 \end{equation}
 Noting that $p^{*}>p$, hence $\|f\|_{L(p^{*},p^{*})}\leq C\|f\|_{L(p^{*},p)}$, (\ref{SLI'}) yields  (\ref{SI}) with $\frac{1}{p^{*}}=\frac{1}{p}-\frac{1}{\sigma}$ and $q< p\leq \sigma$. Using Theorem \ref{kz}, we get (\ref{SI}) for every $q_{0}<p\leq \sigma$ where $q_{0}=\inf\left\lbrace q\in [1,\infty[; (P_{q}) \textrm{ holds }\right\rbrace.$ If $q_{0}=1$, we allow $p=1$.
\begin{rem} 1- As we mentioned in the introduction, a Lie group of polynomial growth satisfies $(D)$, $(P_1)$. Moreover, for $n\in[d,D]$ we have $\mu(B)\geq cr^{n}$ for any ball $B$ of radius $r>0$ --$d$ being the local dimension and $D$ the dimension at infinity--. Therefore this subsection applies on such a group.\\
2-It has been proven \cite{saloff3} that under $(D)$, $(P^{\prime\prime}_{q})$ and $(\ref{GM})$ with $\sigma> q$, the Sobolev inequality (\ref{SI}) holds for all $q\leq p<\sigma$. Since $(D)$ and $(P_{q})$ yield $(P^{\prime\prime}_{q})$, we recover this result under our hypotheses. Besides, we are able to treat the limiting case $p=\sigma$.
 \end{rem}
\section{Proof of Theorem \ref{G} and Theorem \ref{M}} 
The main tool to prove these two theorems is the following two characterizations of the $K$-functional of real interpolation for the homogeneous Sobolev norm.
\begin{thm}\label{EKH}(\cite{badr1}) Let $M$ be a complete Riemannian manifold satisfying $(D)$ and $(P_{q})$ for some $1\leq q<\infty$. Consider the $K$-functional of real interpolation  for the spaces $\overset{.}{W_{q}^{1}}$ and $\overset{.}{W_{\infty}^{1}}$ defined as 
$$
K(F,t,\overset{.}{W_{q}^{1}},\overset{.}{W_{\infty}^{1}})=\inf_{ \substack{f=h+g\\ h\in \overset{.}{E_{q}^{1}},\,g\in\overset{.}{E_{\infty}^{1}}}}\left(\|\,|\nabla h|\,\|_{q}+t\|\,|\nabla g|\,\|_{\infty}\right) 
$$  
where $ f\in \overset{.}{E_{q}^{1}}+\overset{.}{E_{\infty}^{1}}$ with $F=\overline{f}$.
\\
Then
\begin{itemize}
\item[1.] there exists $C_{1}$ such that for every $F \in\overset{.}{W_{q}^{1}}+\overset{.}{W_{\infty}^{1}}$ and $t>0$
$$ K(F,t^\frac{1}{q},\overset{.}{W_{q}^{1}},\overset{.}{W_{\infty}^{1}})\geq C_{1}t^{\frac{1}{q}}|\nabla f|^{q**\frac{1}{q}}(t)\textrm{ where } f\in\overset{.}{E_{q}^{1}}+\overset{.}{E_{\infty}^{1}} \textrm{ with } F=\overline{f};
$$
\item[2.] for $ q\leq p<\infty$, there exists $C_{2}$ such that for every $F\in \overset{.}{W_{p}^{1}}$  and $t>0$
$$ 
K(F,t^{\frac{1}{q}},\overset{.}{W_{q}^{1}},\overset{.}{W_{\infty}^{1}})\leq C_{2} t^{\frac{1}{q}}|\nabla f|^{q**\frac{1}{q}}(t) \textrm{ where } f\in \overset{.}{E_{p}^{1}} \textrm{ such that } F=\overline{f}.
$$
\end{itemize}
\end{thm}
\begin{thm}\label{EK'H} Let $M$ be as in Theorem \ref{EKH}. For $f\in W_{q}^{1}+W_{\infty}^{1}$, consider the functional of interpolation $K'$ defined as follows:
$$
K'(f,t)=K'(f,t,\overset{.}{W_{q}^{1}},\overset{.}{W_{\infty}^{1}})=\inf_{\substack{ f=h+g\\ h\in {W_{q}^{1}},\,g\in{W_{\infty}^{1}}}}\left(\|\,|\nabla h|\,\|_{q}+t\|\,|\nabla g|\,\|_{\infty}\right).
$$  
Let $f\in W_{q}^{1}+W_{\infty}^{1}$ such that $f^{*}(\infty)=0$ and $|\nabla f|^{*}(\infty)=0$. We have
\begin{equation}\label{K'}
K'(f,t^{\frac{1}{q}})\sim t^{\frac{1}{q}}(|\nabla f|^{q**})^{\frac{1}{q}}(t) 
\end{equation}
where the implicit constants do not depend on $f$ and $t$.
Consequently for such $f$'s,
$$ K'(f,t^{\frac{1}{q}})\sim K(\overline{f},t^{\frac{1}{q}},\overset{.}{W_{q}^{1}},\overset{.}{W_{\infty}^{1}}).
$$
\end{thm}
\begin{proof}
Obviously 
$$ t^{\frac{1}{q}}(|\nabla f|^{q**})^{\frac{1}{q}}(t)\leq K(\overline{f},t^{\frac{1}{q}},\overset{.}{W_{q}^{1}},\overset{.}{W_{\infty}^{1}})\leq K'(f,t^{\frac{1}{q}})$$
for all $f\in W_{q}^{1}+W_{\infty}^{1}$.
For the converse estimation, we distinguish three cases:
\begin{itemize}
\item[1.] Let $f\in C_{0}^{\infty}$. For $t>0$, we consider the Calder\'on-Zygmund decomposition given by Proposition 5.5 in \cite{badr1} with $\alpha(t)=\left(\mathcal{M}(|\nabla f|^{q})\right)^{*\frac{1}{q}}(t)\sim(|\nabla f|^{q**})^{\frac{1}{q}}(t)$. We can write then $f=b+g$ with $\|\,|\nabla b|\,\|_{q}\leq C\alpha(t)t^{\frac{1}{q}}$ and $g$ Lipschitz with $\|\,|\nabla g|\,\|_{\infty}\leq C\alpha(t)$ (see also the proof of  Theorem 1.4 in \cite{badr1}). Notice that since $f\in C_{0}^{\infty}$ one has in addition $b\in L_{q}$ and $g\in L_{\infty}$. Consequently, $b\in  W_{q}^{1}$ and $g$ is in $W^{1}_{\infty}$. Therefore, we get (\ref{K'}).
\item[2.] Let $f\in W_{q}^{1}$. There exists a sequence $(f_{n})_{n}$ such that for all $n$, $f_{n}\in C_{0}^{\infty}$ and $\|f-f_{n}\|_{W_{q}^{1}}\rightarrow 0$. Since $|\nabla f_{n}|^{q}\rightarrow |\nabla f|^{q}$ in $L_{1}$, it follows that $|\nabla f_{n}|^{q**}(t)\rightarrow |\nabla f|^{q**}(t)$ for all $t>0$. We have seen in item 1. that for every $n$ there is $g_{n}\in W_{\infty}^{1}$ such that $\|\,|\nabla (f_{n}-g_{n})|\,\|_{q}+t^{\frac{1}{q}}\|\,|\nabla g_{n}|\,\|_{\infty}\leq Ct^{\frac{1}{q}}(|\nabla f_{n}|^{q**})^{\frac{1}{q}}(t)$. Then
 \begin{align*}
\|\,|\nabla(f-g_{n})|\,\|_{q}+t^{\frac{1}{q}}\|\,|\nabla g_{n}|\,\|_{\infty}&\leq \|\,|\nabla(f-f_{n})|\,\|_{q}+\left( \|\,|\nabla (f_{n}-g_{n})|\,\|_{q}+t^{\frac{1}{q}}\|\,|\nabla g_{n}|\,\|_{\infty}\right)
\\
&\leq \epsilon_{n} +Ct^{\frac{1}{q}}(|\nabla f_{n}|^{q**})^{\frac{1}{q}}(t)
\end{align*}
where $\epsilon_{n}\rightarrow 0$ when $n\rightarrow \infty$. We let $n\rightarrow \infty$ to obtain (\ref{K'}).
\item[3.] Let $f\in W_{q}^{1}+W_{\infty}^{1}$ such that $f^{*}(\infty)=0$ and $|\nabla f|^{*}(\infty)=0$. Fix $t>0$ and $p_{0}\in M$. Consider $\varphi\in C_{0}^{\infty}(\mathbb{R})$ satisfying $\varphi\geq 0$, $\varphi(\alpha)=1$ if $\alpha<1$ and $\varphi(\alpha)=0$ if $\alpha>2$. Then put $f_{n}(x)=f(x) \varphi(\frac{d(x,p_{0})}{n})$. Elementary calculations establish that $f_{n}$ lies in $W_{q}^{1}$, hence $K'(f
_{n},t^{\frac{1}{q}})\leq Ct^{\frac{1}{q}}|\nabla f_{n}|^{q**\frac{1}{q}}(t)$. It is shown in \cite{badr1} that
\begin{equation*} K(f,t^{\frac{1}{q}}, W_{q}^{1}, W_{\infty}^{1})\sim \left(\int_{0}^{t}|f|^{q*}(s)ds\right)^{\frac{1}{q}}+\left(\int_{0}^{t}|\nabla f|^{q*}(s)ds\right)^{\frac{1}{q}}.
 \end{equation*}
 All these ingredients yield
  \begin{align}\label{K''}
   K'(f,t^{\frac{1}{q}})&\leq K'(f-f_{n},t^{\frac{1}{q}})+K'(f_{n},t^{\frac{1}{q}})\nonumber
  \\
  &\leq K(f-f_{n},t, W_{q}^{1},W_{\infty}^{1})+K'(f_{n},t)\nonumber
  \\
  &\leq C\left(\int_{0}^{t}|f-f_{n}|^{q*}(s)ds\right)^{\frac{1}{q}}+C\left(\int_{0}^{t}|\nabla f-\nabla f_{n}|^{q*}(s)ds\right)^{\frac{1}{q}}\nonumber
  \\
  &+C\left(\int_{0}^{t}|\nabla f_{n}|^{q*}(s)ds\right)^{\frac{1}{q}}.
  \end{align}
 Now we invoke the following theorem from \cite{krein} page 67-68 stated there in the Euclidean case. As the proof is the same, we state it in the more general case:
 \begin{thm}Let $M$ be a measured space. Consider a sequence of measurable functions $(\psi_{n})_{n}$ and $g$  on $M$ such that $\mu\{|g|>\lambda\}<\infty$ for all $\lambda>0$ with $|\psi_{n}(x)|\leq |g(x)|$. If $\psi_{n}(x)\rightarrow \psi(x) \; \mu-a.e.$ then $(\psi-\psi_{n})^{*}(t)\rightarrow 0\; \forall t>0$. 
 \end{thm}
 We apply this theorem three times:
 \begin{itemize}
 \item[a.] with $\psi_{n}=|f-f_{n}|^{q}$, $\psi=0$ and $g=2^{q}f^{q}$. Using the Lebesgue dominated convergence theorem  we obtain $\int_{0}^{t}|f-f_{n}|^{q*}(s)ds\rightarrow 0$.
 \item[b.] with $\psi_{n}=|\nabla f-\nabla f_{n}|^{q}$, $\psi=0$ and $g=C(|\nabla f|^{q}+|f|^{q})$, where $C$ only depends on $q$, since 
 $$\nabla f_{n}=\nabla f \ind_{B(p_0,n)}+\left(\frac{1}{n}f \varphi'(\frac{d(x,p_{0})}{n})\nabla (d(x,p_{0}))+\nabla f \varphi(\frac{d(x,p_0)}{n})\right)\ind_{B(p_0,n)^{c}}.
 $$
 So again by the Lebesgue dominated convergence theorem we get $\int_{0}^{t}|\nabla f-\nabla f_{n}|^{q*}(s)ds\rightarrow 0$.
 \item[c.] with $\psi_{n}=|\nabla f_{n}|^{q}$, $\psi=|\nabla f|^{q}$ and $g=C(|\nabla f|^{q}+|f|^{q})$, $C$ only depending on $q$, so we get $\int_{0}^{t}|\nabla f_{n}|^{q*}(s)ds\rightarrow \int_{0}^{t}|\nabla f|^{q*}(s)ds$.
  \end{itemize}
  \end{itemize}
  Passing to the limit in (\ref{K''}) yields $K'(f,t^{\frac{1}{q}})\leq C t^{\frac{1}{q}}|\nabla f|^{q**\frac{1}{q}}(t)$ and finishes the proof.
  \end{proof}
\begin{proof}[Proof of Theorem \ref{G}] Let $t>0$, $f\in W_{q}^{1}+W_{\infty}^{1}$ such that $f^{*}(\infty)=0$ and $|\nabla f|^{*}(\infty)=0$. Observe that 
\begin{equation}\label{f}
|f-P_{t}f|^{q**\frac{1}{q}}(s)\leq Ct^{\frac{1}{2}}|\nabla f|^{q**\frac{1}{q}}(s).
\end{equation}
Before proving (\ref{f}), let us see how to conclude from it the desired symmetization inequality. Indeed, (\ref{f}) yields
\begin{align*}
|f|^{q**\frac{1}{q}}(s)&\leq C[\,|f-P_{t}f|^{q**\frac{1}{q}}+|P_{t}f|^{q**\frac{1}{q}}](s)
\\
&\leq C[ t^{\frac{1}{2}}|\nabla f|^{q**\frac{1}{q}}+t^{\frac{\alpha}{2}}t^{-\frac{\alpha}{2}}|P_{t}f|^{q**\frac{1}{q}}](s)
\\
&\leq C t^{\frac{1}{2}}|\nabla f|^{q**\frac{1}{q}}(s)+ Ct^{\frac{\alpha}{2}}\sup\limits_{t>0}\left(t^{-\frac{\alpha}{2}}|P_{t}f|^{q**\frac{1}{q}}(s)\right)
\\
&\leq Ct^{\frac{1}{2}}|\nabla f|^{q**\frac{1}{q}}(s)+Ct^{\frac{\alpha}{2}}\sup\limits_{t>0}t^{-\frac{\alpha}{2}}\|P_{t}f\|_{\infty}
\\
&=C t^{\frac{1}{2}}|\nabla f|^{q**\frac{1}{q}}(s)+Ct^{\frac{\alpha}{2}}\|f\|_{B_{\infty,\infty}^{\alpha}}.
\end{align*}
Therefore we get
\begin{align*}
|f|^{q**\frac{1}{q}}(s)&\leq C\inf\limits_{t>0}\left(t^{\frac{1}{2}}|\nabla f|^{q**\frac{1}{q}}(s)+t^{\frac{\alpha}{2}}\|f\|_{B_{\infty,\infty}^{\alpha}}\right)
\\
&\leq C |\nabla f|^{q**\frac{|\alpha|}{q(1+|\alpha|)}}(s)\|f\|_{B_{\infty,\infty}^{\alpha}}^{\frac{1}{1+|\alpha|}}.
\end{align*}
It remains to prove (\ref{f}). The main tool will be the pseudo-Poincar\'{e} inequalities $(P'_{q})$, $(P'_{\infty})$ and Theorem \ref{EK'H}.

Let $f\in W_{q}^{1}+W_{\infty}^{1}$ such that $f^{*}(\infty)=0$ and $|\nabla f|^{*}(\infty)=0$. Assume that $ f= h+g$ with $ h\in W^{1}_{q},\,g  \in
W^{1}_{\infty}$. We write
$$
f-P_{t}f=(h-P_{t}h)+(g-P_{t}g).
$$
Let $s>0$. The pseudo-Poincar\'e inequalities $(P'_{q})$ and $(P'_{\infty})$ yield
$$
 \|h-P_{t}h\|_{q}+s^{\frac{1}{q}}\|g-P_{t}g\|_{\infty}\leq Ct^{\frac{1}{2}}(\|\,|\nabla h|\,\|_{q}+s^{\frac{1}{q}}\|\,|\nabla g|\,\|_{\infty}).
 $$
 Since $$
K(f,s^{\frac{1}{q}},L_{q},L_{\infty})\sim\Bigl(\int_{0}^{s}(f^{*}(u))^{q}du\Bigr)^{\frac{1}{q}}=s^{\frac{1}{q}}|f|^{q**\frac{1}{q}}(s)
$$
we obtain
\begin{align*}
s^{\frac{1}{q}}|f-P_{t}f|^{q**\frac{1}{q}}(s)&\sim \inf_{\substack{f-P_{t}f=h'+g'\\h'\in L_{q},\, g'\in L_{\infty}}}(\|h'\|_{q}+s^{\frac{1}{q}}\|g'\|_{\infty})
\\
&\leq \inf\limits_{\substack {f=h+g\\ h\in W_{q}^{1},\,g\in W_{\infty}^{1}}}(\|h-P_{t}h\|_{q}+s^{\frac{1}{q}}\|g-P_{t}g\|_{\infty})
\\
&\leq C t^{\frac{1}{2}} \inf\limits_{\substack {f=h+g\\ h\in W_{q}^{1},\, g\in W_{\infty}^{1}}}(\|\,|\nabla h|\,\|_{q}+s^{\frac{1}{q}}\|\,|\nabla g|\,\|_{\infty})
\\
&=Ct^{\frac{1}{2}}K'(f,s^{\frac{1}{q}}).
\end{align*}
Applying Theorem \ref{EK'H}, we obtain the desired inequality (\ref{f}).
 \end{proof}
 \begin{proof}[Proof of Theorem \ref{M}] The proof of this theorem is similar to that of Theorem \ref{G}. Here the key ingredients are the pseudo-Poincar\'{e} inequality for averages $(P^{\prime\prime}_{q})$ that holds for all $f\in \dot{E}_{q}^{1}$. This pseudo-Poincar\'e inequality follows from $(D)$ and the Poincar\'{e} inequality $(P_{q})$. We also make use of Theorem \ref{EK'H}.
 \end{proof}
 \section{Another symmetrization inequality}
 In this section we prove another symmetrization inequality which had been used in \cite{martin} to prove  Gagliardo-Nirenberg inequalities with a Triebel-Lizorkin condition.
 \begin{thm} Let $M$ be a complete non-compact Riemannian manifold satisfying $(D)$ and $(P_{q})$  for some $1\leq q<\infty$. Moreover, assume that $M$ satisfies the pseudo-Poincar\'{e} inequalities $(P'_{q})$ and $(P'_{\infty})$. Consider $\alpha<0$. Then there is $C>0$ such that for every $f\in W_{q}^{1}+W_{\infty}^{1}$ with $f^{*}(\infty)=0$, $|\nabla f|^{*}(\infty)=0$ and satisfying $(\sup\limits_{t>0} t^{-\frac{\alpha}{2}}|P_{t}f(.)|)\in L_{q}+L_{\infty}$:
\begin{equation}\label{P}
|f|^{q**\frac{1}{q}}(s)\leq C |\nabla f|^{q**\frac{|\alpha|}{q(1+|\alpha|)}}(s)\left[\left(\sup\limits_{t>0} t^{-\frac{\alpha}{2}}|P_{t}f(.)|\right)^{q**\frac{1}{q}}(s)\right]^{\frac{1}{1+|\alpha|}},\; s>0.
\end{equation}
\end{thm}
\begin{proof}
 From
$$
|f|^{q} \leq 2^{q-1}\left(|f-P_{t}f|^{q}+t^{\frac{\alpha q}{2}}\sup\limits_{t>0}t^{-\frac{\alpha q}{2}}|P_{t}f|^{q}\right)
 $$
 we obtain
 \begin{align*}
 |f|^{q**\frac{1}{q}}(s) &\leq C\left(|f-P_{t}f|^{q**\frac{1}{q}}(s)+t^{\frac{\alpha}{2}}\left(\sup\limits_{t>0}t^{-\frac{\alpha}{2}}|P_{t}f|\right)^{q**\frac{1}{q}}(s)\right)
 \\
 &\leq C\left(t^{\frac{1}{2}}|\nabla f|^{q**\frac{1}{q}}(s)+t^{\frac{\alpha}{2}}\left(\sup\limits_{t>0}t^{-\frac{\alpha}{2}}|P_{t}f|\right)^{q**\frac{1}{q}}(s)\right).
 \end{align*}
 It follows that
 \begin{align*}
 |f|^{q**\frac{1}{q}}(s) &\leq C \inf\limits_{t>0}\left(t^{\frac{1}{2}}|\nabla f|^{q**\frac{1}{q}}(s)+t^{\frac{\alpha}{2}}\left(\sup\limits_{t>0}t^{-\frac{\alpha}{2}}|P_{t}f|\right)^{q**\frac{1}{q}}(s)\right)
 \\
 &\leq C |\nabla f|^{q**\frac{|\alpha|}{q(1+|\alpha|)}}(s)\left(\sup\limits_{t>0}t^{-\frac{\alpha}{2}}|P_{t}f|\right)^{q**\frac{1}{q(1+|\alpha|)}}(s).
 \end{align*}
 \end{proof}
\section{Proof of Theorem \ref{S}} 
 \begin{proof}
 Let $f\in C^{\infty}_{0}(M)$. Since $p+1\geq 2$, integrating by parts, we get
 $$
 \|\,|\nabla f|\,\|_{p+1}^{p+1}=-\int_{M} div(|\nabla f|^{p-1}\nabla f)fd\mu.
 $$
 Moreover we have $div(|\nabla f|^{p-1}\nabla f)\leq C |\nabla f|^{p-1}|\nabla^{2}f|$. Then
 $$
 \|\,|\nabla f|\,\|_{p+1}^{p+1}\leq C\int_{M}|\nabla f|^{p-1}|\nabla^{2}f||f|d\mu.
 $$
 Let $I=\int_{M}|\nabla f|^{p-1}|\nabla^{2}f||f|d\mu$. Then 
 \begin{align*}
 I&=\int_{0}^{\infty}(|\nabla f|^{p-1}|\nabla^{2}f||f|)^{*}(s)ds
 \\
 &=\int_{0}^{\infty}(|\nabla f|^{\frac{p-2}{2}+\frac{p}{2}}|\nabla^{2}f||f|)^{*}(s)ds
 \\
 &\leq \int_{0}^{\infty}(|\nabla f|^{\frac{p}{2}})^{*}(s)|f|^{q*\frac{1}{q}}(s)(|\nabla f|^{\frac{p-2}{2}}|\nabla^{2}f|)^{*}(s)ds
 \\
 &=\int_{0}^{\infty}|\nabla f|^{q*\frac{p}{2q}}(s)|f|^{q*\frac{1}{q}}(s)(|\nabla f|^{\frac{p-2}{2}}|\nabla^{2}f|)^{*}(s)ds
 \\
 &\leq \int_{0}^{\infty}|\nabla f|^{q**\frac{p}{2q}}(s)|f|^{q**\frac{1}{q}}(s)(|\nabla f|^{\frac{p-2}{2}}|\nabla^{2}f|)^{*}(s)ds.
 \end{align*}
 Thanks to Theorem \ref{G}, we have
 \begin{align*}
 I&\leq \|f\|_{B_{\infty,\infty}^{-1}}^{\frac{1}{2}}\int_{0}^{\infty}|\nabla f|^{q**\frac{p+1}{2q}}(s)( |\nabla f|^{\frac{p-2}{2}}|\nabla^{2}f|^{*})(s)ds
 \\
 &\leq \|f\|_{B_{\infty,\infty}^{-1}}^{\frac{1}{2}}\left(\int_{0}^{\infty}|\nabla f|^{q**\frac{p+1}{q}}(s)ds\right)^{\frac{1}{2}}\left(\int_{0}^{\infty}\left(( |\nabla f|^{\frac{p-2}{2}}|\nabla^{2}f|)^{*}(s)\right)^{2}ds\right)^{\frac{1}{2}}
 \\
 &\leq  \|f\|_{B_{\infty,\infty}^{-1}}^{\frac{1}{2}}\left(\int_{M}|\nabla f|^{p+1}d\mu\right)^{\frac{1}{2}}\left(\int_{M} |\nabla f|^{p-2}|\nabla^{2}f|^{2}d\mu\right)^{\frac{1}{2}}
 \end{align*}
  which finishes the proof.
 \end{proof}
 \begin{rem} Let $M$ be a complete Riemannian manifold satisfying $(D)$ and $(P_{q})$ for some $1\leq q<\infty$. Then Theorem \ref{S} holds replacing the Besov norm $B_{\infty,\infty}^{-1}$ by the Morrey norm $M_{\infty}^{-1}$. This can be proved using Theorem \ref{M}.
 \end{rem}
   
 \bibliographystyle{plain}
\bibliography{latex2}
 \end{document}